\documentclass[11pt]{article} 

\usepackage[utf8]{inputenc} 
\usepackage{url,color}

\usepackage[colorlinks=true,urlcolor=blue,linkcolor=blue,citecolor=magenta]{hyperref}
\usepackage{geometry} 
\usepackage{amsmath}
\usepackage{amssymb,amsthm,nicefrac,paralist}
\usepackage{graphicx} 
\usepackage{mathtools} 
\usepackage{afterpage}
   
\usepackage{booktabs} 
\usepackage{array} 
\usepackage{paralist} 
\usepackage{verbatim}  

\usepackage{tikz}
\usepackage{xcolor}
\usetikzlibrary{shapes,shadows,arrows,patterns,trees,automata,positioning,snakes,fadings,backgrounds,intersections,calc}
\usepackage{tkz-euclide}
\usepackage{ifthen}
\tikzstyle{edge}=[line width=1mm, color=blue] 
\tikzstyle{face}=[color=black,fill=white,opacity=0.8,line width=0.7mm]
\tikzstyle{line}=[black!50!white, opacity=1, line width=0.45mm] 
\tikzstyle{label} = [color=black, font=\normalsize, opacity=1]
\tikzstyle{vertex} = [outer color=black,draw, color=black,line width=0.03mm, inner color=black!55!white, circle,inner sep=0.7mm,minimum size=0.9mm]

\theoremstyle{plain}
\newtheorem{theorem}{Theorem} 

\newtheorem{lemma}[theorem]{Lemma}

\newtheorem*{theorem*}{Theorem}

\theoremstyle{definition}

\newcommand\R{\mathbb{R}}
\newcommand\NC{\mathrm{NC}}
\newcommand\NCC{\mathrm{NC}^c}
\newcommand\CG{\mathrm{G}}

\title{Simple polytopes without small separators, {II}:\\ {Thurston}'s bound\footnote{The first author was funded by DFG through the \emph{Berlin Mathematical School}. Research by the second author was supported by the DFG Collaborative Research Center TRR~109 ``Discretization in Geometry and Dynamics.''}}
\author{Lauri Loiskekoski\\ 
Institut f\"ur Mathematik, FU Berlin\\Arnimallee 2\\14195 Berlin, Germany\\
\url{lauri.loiskekoski@fu-berlin.de}
\and  
G\"unter M. Ziegler\\
Institut f\"ur Mathematik, FU Berlin\\Arnimallee 2\\14195 Berlin, Germany\\
\url{ziegler@math.fu-berlin.de}}
\date{August 22, 2017}

\begin{document}
\maketitle

\begin{abstract}
We show that there are simple $4$-dimensional polytopes with $n$ vertices 
such that all separators of the graph have size at least $\Omega(n/\log n)$. 
This establishes a strong form of a claim by Thurston, for which
the construction and proof had been lost.  

We construct the polytopes by cutting off the vertices and then the edges of 
a particular type of neighborly cubical polytopes.
The graphs of simple polytopes thus obtained are $4$-regular;
they contain $3$-regular ``cube-connected cycle graphs'' as minors of spanning subgraphs.
\end{abstract}

\section{Introduction}

For an arbitrary fixed constant~$c$ with $0<c<\nicefrac12$
(where traditionally $c=\nicefrac13$ is used)
and a simple graph with vertex set $V$,
a \emph{separator} with \emph{separation constant} $c$
is a partition of $V$ into sets $A,B,C$ with 
$cn\le|A|\le|B|\le(1-c)n$ such that $C$ separates $A$ from $B$, that is,
there is no edge between a vertex in $A$ and a vertex in~$B$. 
The \emph{size} of the separator is defined as the size of the set $C$.
 
The Lipton--Tarjan planar separator theorem from 1979 \cite{LiTa}
states that each planar graph on $n$ vertices---and 
in particular the graph of any $3$-dimensional polytope---has
a separator of size $|C|=O(\sqrt n)$.
No such general result can exist for $d$-dimensional polytopes,
as their graphs may be complete. However, the situation for simple 
$d$-polytopes, whose graphs are $d$-regular, is interesting:
Kalai had conjectured in \cite[Conj.~20.2.12]{kalai04:_polyt}
that for any fixed $d\ge2$,
the graphs of all simple $d$-polytopes on $n$ vertices have separators of size
\[
  O\left( n^{1-\frac1{d-1}}\right).
\]

In our previous paper \cite{Z150} we have shown that
cutting off the vertices, and then the edges, from \emph{any} neighborly cubical $4$-polytope
(a cubical $4$-dimensional polytope with the graph of an $m$-cube)
yields a simple $4$-dimensional polytope with $n=\Theta(m2^m)$ vertices
such that any vertex separator has size at least $\Omega(n/\log^{3/2} n)$.

This disproved Kalai's conjecture,  
but it fell short of establishing a claim by Thurston, documented by 
Kalai in \cite[p.~460]{kalai04:_polyt}, who refers to \cite{MillerTengThurstonVavasis}.
Indeed this does not appear there, 
but Gary Miller told us (personal communication) that
in the context of the work towards~\cite{MillerTengThurstonVavasis}, Thurston
had claimed that there are simplicial $3$-spheres on $n$ tetrahedra
for which all separators of the dual graph contain at least $\Omega(n/\log n)$ vertices.
Miller reports that for this 
“\emph{Thurston gave an embedding of the cube-connected cycle graph in $\R^3$ as linear tets}
[tetrahedra]”. Thurston's construction for this seems to be lost.

In this paper, we prove a stronger version of this:
There is not only an embedding of tetrahedra in $\R^3$, or 
a triangulated simplicial $3$-sphere with $n$ facets for which 
every separator of the dual graph has size at least $\Omega(n/\log n)$, 
but indeed there is a simplicial convex $4$-polytope with this property.
Our main theorem describes the dual polytope, which is simple:

\begin{theorem}
	There is a family of simple convex $4$-dimensional polytopes
	$\NCC_4(m)''$ with $n=\Theta(m2^m)$ vertices
	for which the smallest separators of the graph have size $\Theta(2^m)=\Theta(n/\log n)$. 
\end{theorem}

We will obtain such polytopes by
cutting off first all the vertices, and then the edges, from a \emph{specific} sequence of
neighborly cubical $4$-polytopes.
In Section~\ref{sec:construction} we describe the construction of the
polytopes $\NCC_4(m)$, $\NCC_4(m)'$, and finally $\NCC_4(m)''$, 
and derive the combinatorial properties, and in particular 
collect all information that will be needed on their graphs, which we denote by $\CG_m$, $\CG_m'$, and $\CG_m''$, respectively. 
The graphs $\CG_m'$ contain  as spanning subgraphs 
the cube-connected cycle graphs of Preparata \& Vuillemin \cite{PreparataVuillemin1981} that
Thurston presumably referred to.

Then in Section~\ref{sec:expansion} we review Sinclair's canonical paths method,
which allows one to show that a graph has large edge expansion,
which in turn implies that all separators are large.

Finally, in Section~\ref{sec:paths} we specify canonical paths in $\CG_m''$
and show that not too many paths are routed through any edge of any of the four different types of edges in $\CG''$,
which we call ``long,'' ``medium,'' ``extra,'' and ``short''. This completes the proof.

\section{Neighborly cubical 4-polytopes and their vertex and edge truncations}\label{sec:construction}

\subsection{Neighborly cubical polytopes}

A \emph{neighborly cubical $d$-polytope} 
is a $d$-dimensional polytope $\NC_d(m)$ that is cubical (that is, its facets are 
combinatorially equivalent to the $(d-1)$-cube $[0,1]^{d-1}$) and whose graph
is isomorphic to the graph of an $m$-cube, for some $m\ge d\ge4$.
(Except for the last section, we will in this paper mainly refer to the case $d=4$.) 

The boundary complex of any such polytope is combinatorially equivalent to a subcomplex
of the $m$-cube $[0,1]^m$, so its non-empty faces can be identified with vectors in~$\{0,1,*\}^m$,
or equivalently in $\{-,+,*\}^m$.
 
The definition of neighborly cubical $4$-polytopes immediately yields $f_0$ and $f_1$.
Using the Euler equation and double-counting, which yields $3f_3=f_2$, 
we conclude that every neighborly cubical $4$-polytope has the $f$-vector
\[ 
    f(\NC_4(m)) = (f_0,f_1,f_2,f_3) = 2^{m-2} (4, 2m, 3m-6, m-2). 
\]

\subsection{The cyclic neighborly cubical polytopes}

Joswig \& Ziegler \cite{Z62} established the existence of neighborly cubical $d$-polytopes
on $2^m$ vertices, for all $m\ge d\ge4$,  
by constructing one specific example of a neighborly cubical $d$-polytope $\NCC_d(m)$,
which we here call a \emph{cyclic neighborly cubical polytope} $\NCC_d(m)$.
This name is chosen to reflect that -- according to the analysis of Sanyal \& Ziegler \cite{Z102}
-- all vertex figures of $\NCC_d(m)$ are combinatorially isomorphic to a pyramid over a 
triangulation of the cyclic polytope $C_{d-2}(m-1)$.

The polytopes $\NCC_d(m)$ were originally obtained by a somewhat subtle
linear algebra/matrix theory construction.
However, the ``Cubical Gale Evenness Criterion'' \cite[Thm.~18]{Z62} provides 
a complete description
of the combinatorics of the polytope $\NCC_4(m)$.

The graph $\CG_m$ of $\NCC_d(m)$ is isomorphic to the graph of the $m$-cube.
Thus we label each vertex with an entry in $\{{-},{+}\}^m$. 
The edges are labeled by vectors in $\{{-},{+},*\}^m$ with exactly one $*$-entry. 
The coordinate $i\in[m]=\{1,\dots,m\}$ with the  $*$-entry will be called 
the \emph{direction} of the edge. 
There are $2^{m-1}$ edges in each direction.

We now specialize to the case $d=4$. 
The Cubical Gale Evenness Criterion is quite complicated even in this case,
but we will not need the full description.
 
\begin{theorem}
[Part of the Cubical Gale Evenness Criterion, for $d=4$ {\cite[Thm.~18]{Z62}}] 
The facets of the cyclic neighborly polytope $\NCC_4(m)$ are given by vectors in
$\{-,+,*\}^m$ with exactly three $*$s. 

If the first component of a vector in $\{-,+,*\}^m$ is $*$, then it corresponds to a facet 
of $\NCC_4(m)$ if and only if the rest of the vector satisfies the Gale Evenness criterion, 
that is, if between any two non-$*$-entries there is an even number of $*$-entries.
Equivalently, this happens if in the rest of that vector the two $*$-entries are
cyclically adjacent.
\end{theorem}

This theorem shows that at any vertex $v\in\{-,+\}^m$,
the two incident edges in directions $i$ and $i+1 \pmod m$ span a $2$-face, as
$* v_2 v_3 \dots v_{i-1} {*} {*} v_{i+2} \dots v_n\in\{-,+,*\}^m$ 
corresponds to a facet of $\NCC_4(m)$ for $2\le i<m$.
As the $2$-faces correspond to edges of the vertex figure, this implies that the vertices of the vertex figure are cyclically connected in a consistent way independent of which vertex figure we are looking at: There is a Hamilton cycle
$123\dots m$ in the graph of each vertex figure.

Thus the boundary complex of $\NCC_4(m)$ contains a copy of the polyhedral surface
described by Ringel \cite{ringel55:_ueber_probl_wuerf_wuerf} in 1955; this was pointed out
explicitly in Ziegler \cite[Sect.~3]{Z100}.
See also Joswig \& Rörig \cite{joswig:_neigh}.

\subsection{Truncating the vertices}

The polytope $\NCC_4(m)'$ is obtained by
truncating all the vertices of the cyclic neighborly cubical polytope $\NCC_4(m)$.  

The resulting polytope has
\begin{compactitem}
	\item $(m-2)2^{m-2}$ facets that arise from the ``old'' facets (of $\NCC_4(m)$);
		which are vertex-truncated cubes (and hence simple $3$-polytopes),
		and 
	\item $2^m$ ``new'' facets, which are simplicial $3$-polytopes with $m$ vertices each.
\end{compactitem} 
Moreover, the vertices of each ``new'' facet are naturally labeled by $1,2,\dots,m$, 
where the label is given by the direction of the edge of $\CG_m$ that the new vertex lies on.
 
Thus the graph $\CG_m'$ of $\NCC_4(m)'$ -- illustrated in the middle of Figure~\ref{figure:graphs} --
has $m2^m$ vertices, which may be labeled by $(v,i)\in \{0,1\}^m \times [m]$,
where $v\in\{0,1\}^m$ denotes the vertex of $\NCC_4(m)$ which has been truncated,
and $i\in[m]$ is the direction of the edge of $\NCC_4(m)$ which has been cut.

\afterpage{%
\begin{figure}

%
%
%
%
%
%
%
	
\begin{tikzpicture}[scale=2.5,z={(-3mm,-1.5mm)},line join=round,line cap=round]



\clip (-1.5,-.6) rectangle (3.5,1);

\def\xs{0.2}
\def\ys{0.1}
\def\zs{0.1}


\node (A) at (1+\xs,0+\ys,-1+\zs) {};
\node (B) at (-0.5+\xs,0+\ys,-1+\zs) {};
\node (C) at (0.5+\xs,-1+\ys,0+\zs) {};
\node (D) at (1+\xs,-0.02+\ys,1+\zs) {};
\node (E) at (0.5+\xs,1+\ys,0+\zs) {};


\node (a) at ($(0,0,0)!5!(A)$) {};
\node (b) at ($(0,0,0)!10!(B)$) {};
\node (c) at ($(0,0,0)!5!(C)$) {};
\node (d) at ($(0,0,0)!5!(D)$) {};
\node (e) at ($(0,0,0)!5!(E)$) {};


\node (Ab) at ($(A.center)!0.17!(B.center)$) {};
\node (Ac) at ($(A.center)!0.22!(C.center)$) {};
\node (Ad) at ($(A.center)!0.15!(D.center)$) {};
\node (Ae) at ($(A.center)!0.22!(E.center)$) {};

\node (Ba) at ($(B.center)!0.06!(A.center)$) {};
\node (Bc) at ($(B.center)!0.23!(C.center)$) {};
\node (Bd) at ($(B.center)!0.27!(D.center)$) {};
\node (Be) at ($(B.center)!0.23!(E.center)$) {};

\node (Ca) at ($(C.center)!0.29!(A.center)$) {};
\node (Cb) at ($(C.center)!0.29!(B.center)$) {};
\node (Cd) at ($(C.center)!0.29!(D.center)$) {};

\node (Da) at ($(D.center)!0.23!(A.center)$) {};
\node (Db) at ($(D.center)!0.25!(B.center)$) {};
\node (Dc) at ($(D.center)!0.15!(C.center)$) {};
\node (De) at ($(D.center)!0.23!(E.center)$) {};

\node (Ea) at ($(E.center)!0.29!(A.center)$) {};	
\node (Eb) at ($(E.center)!0.29!(B.center)$) {};	
\node (Ed) at ($(E.center)!0.29!(D.center)$) {};



\draw[line] (0,0,0) -- (a.center);
\draw[line] (0,0,0) -- (b.center);
\draw[line] (0,0,0) -- (c.center);
\draw[line] (0,0,0) -- (d.center);
\draw[line] (0,0,0) -- (e.center);

\node[vertex] at (0,0,0) {};

\node at (-1.2,0.5,0) {\large$\CG_5$:};

\end{tikzpicture}

\bigskip


%
%
%
%
%
%
%
	
\begin{tikzpicture}[scale=2.5,z={(-3mm,-1.5mm)},line join=round,line cap=round]



\clip (-1.5,-1.3) rectangle (3.5,1.8);

\def\xs{0.2}
\def\ys{0.1}
\def\zs{0.1}


\node (A) at (1+\xs,0+\ys,-1+\zs) {};
\node (B) at (-0.5+\xs,0+\ys,-1+\zs) {};
\node (C) at (0.5+\xs,-1+\ys,0+\zs) {};
\node (D) at (1+\xs,-0.02+\ys,1+\zs) {};
\node (E) at (0.5+\xs,1+\ys,0+\zs) {};


\node (a) at ($(0,0,0)!5!(A)$) {};
\node (b) at ($(0,0,0)!10!(B)$) {};
\node (c) at ($(0,0,0)!5!(C)$) {};
\node (d) at ($(0,0,0)!5!(D)$) {};
\node (e) at ($(0,0,0)!5!(E)$) {};


\node (Ab) at ($(A.center)!0.17!(B.center)$) {};
\node (Ac) at ($(A.center)!0.22!(C.center)$) {};
\node (Ad) at ($(A.center)!0.15!(D.center)$) {};
\node (Ae) at ($(A.center)!0.22!(E.center)$) {};

\node (Ba) at ($(B.center)!0.06!(A.center)$) {};
\node (Bc) at ($(B.center)!0.23!(C.center)$) {};
\node (Bd) at ($(B.center)!0.27!(D.center)$) {};
\node (Be) at ($(B.center)!0.23!(E.center)$) {};

\node (Ca) at ($(C.center)!0.29!(A.center)$) {};
\node (Cb) at ($(C.center)!0.29!(B.center)$) {};
\node (Cd) at ($(C.center)!0.29!(D.center)$) {};

\node (Da) at ($(D.center)!0.23!(A.center)$) {};
\node (Db) at ($(D.center)!0.25!(B.center)$) {};
\node (Dc) at ($(D.center)!0.15!(C.center)$) {};
\node (De) at ($(D.center)!0.23!(E.center)$) {};

\node (Ea) at ($(E.center)!0.29!(A.center)$) {};	
\node (Eb) at ($(E.center)!0.29!(B.center)$) {};	
\node (Ed) at ($(E.center)!0.29!(D.center)$) {};



%
%
%
%


\draw[face] (A.center) -- (B.center) -- (C.center) -- cycle;
\draw[face] (A.center) -- node[label, below] {\qquad \quad \ medium} (C.center) -- (D.center) -- cycle;
\draw[face] (B.center) -- (C.center) -- (D.center) -- cycle;

\draw[face] (A.center) -- (B.center) -- (E.center) -- cycle;

\draw[edge] (B.center) -- (A.center);

\draw[face] (A.center) -- node[label, above] {\qquad \ extra} (E.center) -- (D.center) -- cycle;
\draw[face] (B.center) -- (E.center) -- (D.center) -- cycle;

\draw[edge] (B.center) -- (E.center) -- (D.center) -- (C.center) -- (A.center);

\draw[line] (A.center) -- (a.center);
\draw[line] (B.center) -- node[label, left] {long edge} (b.center);
\draw[line] (C.center) -- (c.center);
\draw[line] (D.center) -- (d.center);
\draw[line] (E.center) -- (e.center);

\path[every node/.append style=vertex] node at (A) {} node at (B) {} node at (C) {} node at (D) {} node at (E) {};

\draw[line] (D.center) -- (d.center);

\node at (-1.2,0.5,0) {\large$\CG_5'$:};

\end{tikzpicture}

\bigskip


%
%
%
%
%
%
%
	
\begin{tikzpicture}[scale=2.5,z={(-3mm,-1.5mm)},line join=round,line cap=round]



\clip (-1.5,-1.3) rectangle (3.5,1.8);

\def\xs{0.2}
\def\ys{0.1}
\def\zs{0.1}


\node (A) at (1+\xs,0+\ys,-1+\zs) {};
\node (B) at (-0.5+\xs,0+\ys,-1+\zs) {};
\node (C) at (0.5+\xs,-1+\ys,0+\zs) {};
\node (D) at (1+\xs,-0.02+\ys,1+\zs) {};
\node (E) at (0.5+\xs,1+\ys,0+\zs) {};


\node (a) at ($(0,0,0)!5!(A)$) {};
\node (b) at ($(0,0,0)!10!(B)$) {};
\node (c) at ($(0,0,0)!5!(C)$) {};
\node (d) at ($(0,0,0)!5!(D)$) {};
\node (e) at ($(0,0,0)!5!(E)$) {};


\node (Ab) at ($(A.center)!0.17!(B.center)$) {};
\node (Ac) at ($(A.center)!0.22!(C.center)$) {};
\node (Ad) at ($(A.center)!0.15!(D.center)$) {};
\node (Ae) at ($(A.center)!0.22!(E.center)$) {};

\node (Ba) at ($(B.center)!0.06!(A.center)$) {};
\node (Bc) at ($(B.center)!0.23!(C.center)$) {};
\node (Bd) at ($(B.center)!0.27!(D.center)$) {};
\node (Be) at ($(B.center)!0.23!(E.center)$) {};

\node (Ca) at ($(C.center)!0.29!(A.center)$) {};
\node (Cb) at ($(C.center)!0.29!(B.center)$) {};
\node (Cd) at ($(C.center)!0.29!(D.center)$) {};

\node (Da) at ($(D.center)!0.23!(A.center)$) {};
\node (Db) at ($(D.center)!0.25!(B.center)$) {};
\node (Dc) at ($(D.center)!0.15!(C.center)$) {};
\node (De) at ($(D.center)!0.23!(E.center)$) {};

\node (Ea) at ($(E.center)!0.29!(A.center)$) {};	
\node (Eb) at ($(E.center)!0.29!(B.center)$) {};	
\node (Ed) at ($(E.center)!0.29!(D.center)$) {};



%
%
%
%
%


%
%
%
%
%
%
%
%
%
%
%


\draw[line] (Ab.center) -- +($(a)-(A)$);
\draw[line] (Ac.center) -- +($(a)-(A)$);
\draw[line] (Ad.center) -- +($(a)-(A)$);
\draw[line] (Ae.center) -- +($(a)-(A)$);

\draw[line] (Ca.center) -- +($(c)-(C)$);
\draw[line] (Cb.center) -- +($(c)-(C)$);
\draw[line] (Cd.center) -- +($(c)-(C)$);

\draw[face] (Ca.center) -- (Cb.center) -- (Cd.center) -- cycle;

\draw[face] (Ac.center) -- (Ab.center) -- (Ba.center) -- node[label, left] {short} (Bc.center) -- (Cb.center) -- (Ca.center) -- cycle;

\draw[face] (Bd.center) -- (Bc.center) -- (Cb.center) -- (Cd.center) -- (Dc.center) --(Db.center) -- cycle;

\draw[face] (Ae.center) -- (Ab.center) -- (Ba.center) -- (Be.center) -- (Eb.center) -- (Ea.center) -- cycle;

\draw[face] (Ab.center) -- (Ae.center) -- (Ad.center) -- (Ac.center) -- cycle;

\draw[face] (Ad.center) -- (Ac.center) -- node[label, below] {\qquad \quad \ medium} (Ca.center) -- (Cd.center) -- (Dc.center) -- (Da.center) -- cycle;

\draw[edge] (Ba.center) -- (Ab.center);

\node[vertex] at (Ab) {};

\draw[face] (Da.center) -- (De.center) -- (Db.center) -- (Dc.center) -- cycle;

\draw[face] (Ad.center) -- (Ae.center) -- node[label, above] {\qquad \ extra} (Ea.center) -- (Ed.center) -- (De.center) -- (Da.center) -- cycle;

\draw[face] (Ba.center) -- (Bc.center) -- (Bd.center) -- (Be.center) -- cycle;

\draw[face] (Bd.center) -- (Be.center) --  (Eb.center) -- (Ed.center) -- (De.center) -- (Db.center) -- cycle;

\draw[face] (Ea.center) -- (Eb.center) -- (Ed.center) -- cycle;

\draw[edge] (Be.center) -- (Eb.center);
\draw[edge] (Ed.center) -- (De.center);
\draw[edge] (Dc.center) -- (Cd.center);
\draw[edge] (Ca.center) -- (Ac.center);

\node[vertex] at (Ac) {};

\draw[line] (Ba.center) -- node[label, left] {long edges} +($(b)-(B)$);
\node[vertex] at (Ba) {};
\draw[line] (Bc.center) -- +($(b)-(B)$);
\draw[line] (Bd.center) -- +($(b)-(B)$);
\draw[line] (Be.center) -- +($(b)-(B)$);

\draw[line] (Ea.center) -- +($(e)-(E)$);
\draw[line] (Eb.center) -- +($(e)-(E)$);
\draw[line] (Ed.center) -- +($(e)-(E)$);

\path[every node/.append style=vertex] node at (Ad) {} node at (Ae) {} node at (Bc) {} node at (Bd) {} node at (Be) {} node at (Ca) {} node at (Cb) {} node at (Cd) {} node at (Db) {} node at (Dc) {} node at (De) {} node at (Ea) {} node at (Eb) {} node at (Ed) {};

\node at (-1.2,0.5,0) {\large$\CG_5''$:};

\node[vertex] at (Da) {};
\draw[line] (Da.center) -- +($(d)-(D)$);

\draw[line] (Db.center) -- +($(d)-(D)$);
\draw[line] (Dc.center) -- +($(d)-(D)$);
\draw[line] (De.center) -- +($(d)-(D)$);

\end{tikzpicture}

	
\caption{The graphs $\CG_m$, $\CG_m'$, and $\CG_m''$, for $m=5$:
	Local situation at one cluster.}
\label{figure:graphs}
\end{figure}
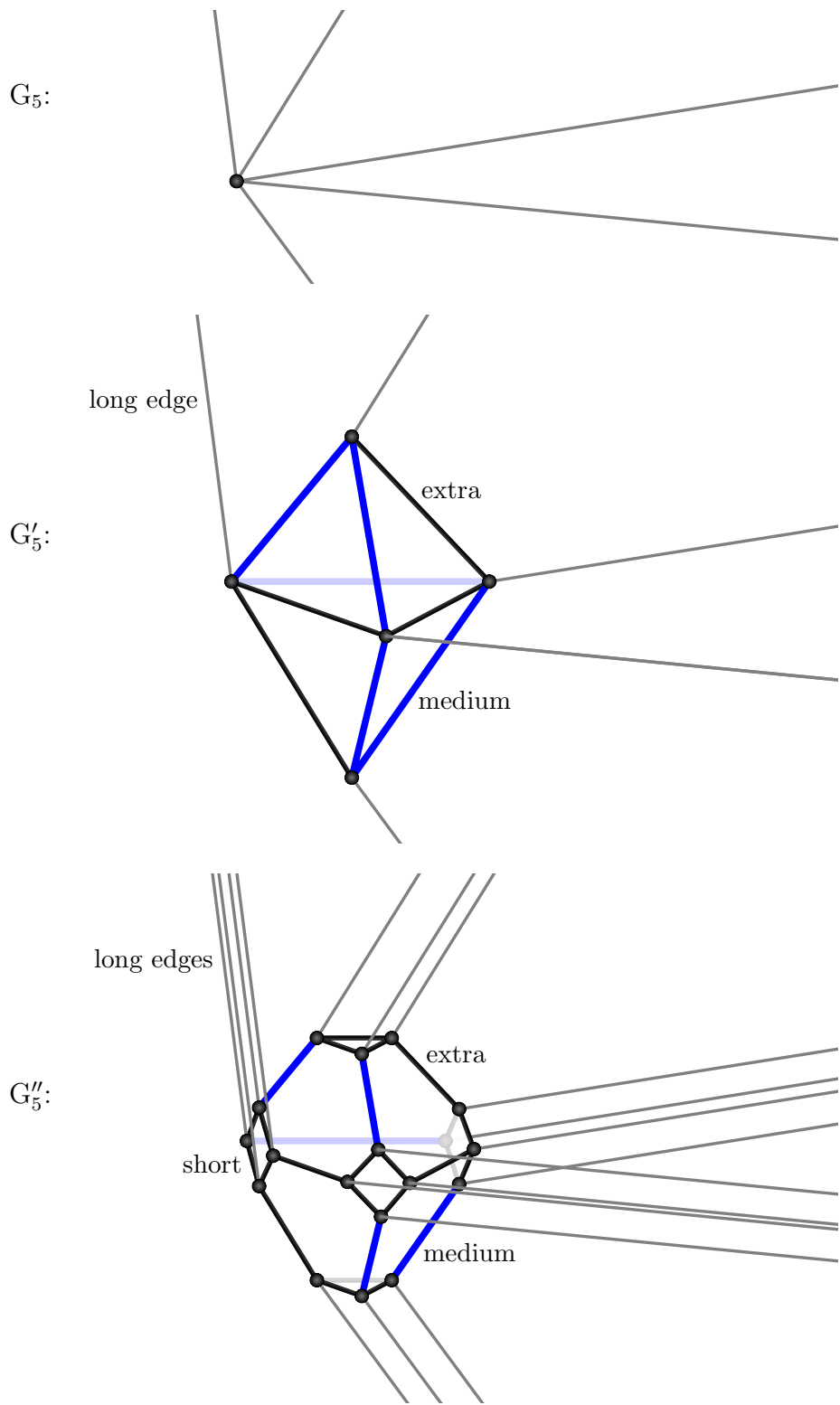	
\clearpage
}

Indeed, there are three types of edges in $\CG_m'$:
First, there are the $m2^{m-1}$ ``long'' edges between vertices $(v,i)$ and $(v',i)$,
where $v$ and $v'$ differ only in coordinate $i$.
Then there are the $m2^m$ ``medium'' edges between vertices $(v,i)$ and $(v,i+1\pmod m)$
on the Hamilton cycle on the graph of the vertex figure at $v$,
which is given by $(1,2,\dots,m)$. The ``long'' and ``medium'' edges together form the cube-connected 
cycle graphs of Preparata \& Vuillemin \cite{PreparataVuillemin1981}. 
Finally there are $(2m-6)2^m$ additional ``extra'' edges between vertices $(v,i)$ and $(v,j)$
with $j-i\neq\pm1\pmod m$, which we will not use in the following, but which together with the
``medium'' edges yield the graphs of the $2^m$ (simplicial) vertex figures. 

In describing $\CG_m'$ we refer to the $2^m$ subgraphs formed by the ``medium'' and ``extra''
edges as \emph{clusters}. These cluster subgraphs are maximal planar graphs on $m$ vertices
and $3m-6$ edges. The clusters are connected by ``long'' edges in a cube-like fashion. 

\subsection{Then cutting off the edges}\label{subsec:NC2}\enlargethispage{3mm}%

The polytope $\NCC_4(m)''$ is then obtained from $\NCC_4(m)'$ by cutting off what
remains of the $m2^{m-1}$ original edges of $\NCC_4(m)$.
The new facets formed by this have the combinatorial type of 
prisms over $k$-gons, with $3\le k\le m-1$.
Here $k$ is the number of facets that contained the edge that was cut off,
or equivalently the degree of its end vertices in the cluster triangulations
(that is, in the vertex figures from the previous step).
The graph of each prism consists of two \emph{cycles of short edges},
which both are $k$-cycles, and of set of $k$ long edges that we refer
to as a \emph{parallel class of long edges}. 

The graph $\CG_m''$ of the resulting polytope $\NCC_4(m)''$ 
 -- illustrated in the bottom part of Figure~\ref{figure:graphs} --
has $2^m$ planar subgraphs that we again refer to as \emph{clusters}, which are connected in an $m$-cube like pattern
by parallel classes of at least $3$ and at most $m-1$ ``long'' edges.
Each cluster in $\CG_m''$ is a $3$-connected $3$-regular planar graph on $6m-12$ vertices,
which consists of $m$ cycles formed by ``short'' edges and
$m$ ``medium'' and $2m-6$ ``extra''  edges.

The cycles of short edges correspond bijectively to the vertices of $\CG_m'$; 
thus they get labels of the form $(v,i)$.
``Medium'' edges connect the subsequent cycles $(v,i)$ and $(v,j)$
with $j=i+1\pmod m$, while the ``extra'' edges connect non-subsequent cycles.

\subsection{Minor relations} 

There are minor relations 
\[
   \CG_m''\rightarrow \CG_m'\rightarrow\CG_m.
\]
Here $\CG_m'$ is obtained from $\CG_m''$ by contraction of the ``short''
edges in $\CG_m''$ and identification of ``long'' edges that become parallel after the contractions.
Similarly, $\CG_m$ is obtained from $\CG_m'$ by contraction of the ``medium'' and ``extra''
edges in $\CG_m'$. 

\subsection{Greater generality}
 
The construction and analysis in this section
can be extended to the large number of combinatorial types 
of neighborly cubical $4$-polytopes described and analyzed 
by Sanyal \& Ziegler \cite{Z102}:
We can find a cycle going through the directions in a
cluster for any labeling of the vertices of the polygon; \cite[Thm.~3.7]{Z102} 
tells how the
triangulation is done. In particular, the combinatorial type of the vertex figure
for a vertex $v \in \NC_4(m)$ is created by pushing
the vertices of the polygon in the specified order
until we hit a plus sign in the vertex label, which corresponds to pulling
and completes the triangulation. Changing the order of the vertices means
our cycle would not be $(1,2,\dots,m)$, but some other (fixed!) permutation of~$[m]$. 

\section{Edge expansion and vertex separators for simple graphs}\label{sec:expansion}

\subsection{Definitions}

Let $G = (V,E)$ be a simple graph. For any set $S \subset V$, the \emph{edge boundary} 
$\delta(S)$ is the set of edges with one end in $S$ and the other in 
$V{\setminus}S$. The \emph{edge expansion} $\mathcal{X}(G)$ of $G$ is then defined by
\begin{align}
 \mathcal{X}(G) := \min \left\{ 
\frac{|\delta(S)|}{|S|}\,:\, S\subset V,\ S\neq \emptyset,\ S\le \frac{|V|}{2} \right\}.
\end{align}
This quantity is crucial for the mixing properties of Markov chains on~$G$.
We consider it because of its close connection to separators: For regular graphs with large
expansion, separators have to be large as well (see below).

For some examples (including ours), the expansion can very efficiently be estimated by 
the method of ``canonical paths'' pioneered by
Sinclair~\cite{Sinclair}. A nice introduction and overview was given by Kaibel~\cite{Kaibel}.
(A number of powerful generalizations of Sinclair's method exist, for example
using random paths or, equivalently, flows. We will not use these.)

\subsection{Sinclair's canonical paths method}\label{subsec:Sinclair}

To estimate the edge expansion, Sinclair tells us to specify,
for each pair of vertices $(s,t) \in V \times V$,
a path from $s$ to $t$.
Let $\phi: E(G) \rightarrow \mathbb{N}$ count for each edge the number of paths 
that use it. We can now bound the edge expansion in terms of 
\[
	\phi_{\max} := \max\{\phi(e) : e\in E(G)\}.
\] 

\begin{lemma}[Sinclair's lemma]
Let $\phi_{\max}$ be defined as in the previous paragraph. 
Then the edge expansion of the graph $G$ satisfies
\[
\mathcal{X}(G) \ge \frac{n}{2\,\phi_{\max}}.
\]
\end{lemma}

\begin{proof}
	Let $\phi(\delta(S))$ be the sum of $\phi(e)$ over all edges in $\delta(S)$,
which is at least the number of paths that use edges in $\delta(S)$. 
Thus we get $\phi(\delta(S)) \ge |S|\,|V{\setminus}S|$. On the other 
hand we clearly have
$\phi(\delta(S)) \le \phi_{\max} \cdot |\delta(S)|$. This means that 
for any $|S|\le \frac{|V|}{2}$,
\[
 \mathcal{X}(G) \ge \frac{|\delta(S)|}{|S|} 
                \ge \frac{\phi(\delta(S))}{\phi_{\max}\,|S|}
                \ge \frac{|S|\,|V{\setminus}S|}{\phi_{\max}\,|S|}
                 =  \frac{|V{\setminus}S|}{\phi_{\max}}
				\ge \frac{n}{2\,\phi_{\max}}.
\]
\vskip-7mm
\end{proof}

\subsection{Relating edge expansion and vertex separators}\label{subsec:expansion}
 
\begin{lemma}
	In a $d$-regular graphs $G$ on $n$ vertices with edge expansion $\mathcal{X}(G)$,
	all separators have size at least
	\[
	\frac cd\mathcal{X}(G)n = \Omega(\mathcal{X}(G)n),
	\]
	where $d$ is the constant from the definition of a separator.
\end{lemma}

\begin{proof}
Let $G$ be a $d$-regular graph with expansion $\mathcal{X}(G)$  
and let $(A,B,C)$ be a separator of $G$, with $|B| \ge|A| \ge cn$.  
By definition there are at least $\mathcal{X}(G)|A|$ edges in the boundary $\delta(A)$. 
The other ends of these edges lie in~$C$, as there are no edges between $A$ and $B$. 
Since $G$ is $d$-regular, $|C|$ has size at least $\mathcal{X}(G)|A|/d \ge \mathcal{X}(G)cn/d = (c/d)\mathcal{X}(G)n$.  
\end{proof}

\section{Canonical paths}\label{sec:paths}

\subsection{Defining the paths}

In order to lower-bound the expansion of $\CG_m''$, 
we now define paths between every pair of vertices of $\NCC_4(m)$. 

Let $v_0$ be a vertex in the cycle of short edges labeled $(v,i)$
and $w_0$ be a vertex in the cycle labeled $(w,j)$. 
Now route the path from $v_0$ to $w_0$ as follows.
For each of the coordinates, taken in cyclic order $i,i+1,\dots,m,1,\dots,i-1$,
perform the following procedure: 
	
\begin{quote}
	\textbf{Procedure P:}
	If the path as constructed up to now ends at a vertex $u_0\in(u,k)$, then
\begin{compactitem}
	\item If $u$ and $w$ differ in coordinate $k$, then take the
	      long edge incident to $u_0$, which leads to the cycle $(u',k)$, where $u$ and $u'$ 
		  differ only in the coordinate $k$. 
		  Then use only short edges and then a medium edge to get to
		  the cycle labeled $(u',k+1)$.
	\item If $u$ and $w$ do not differ in coordinate $k$, then 
		  use only short edges and then one medium edge to get to
		  the cycle labeled $(u,k+1)$.
\end{compactitem}
\end{quote}

This procedure is performed at most $m$ times until the path constructed
reaches a vertex in the cluster labelled $w$, and then at most $m-1$ times
within this cluster until the path constructed reaches a vertex in the cycle labelled $(w,j)$,
and then in one final iteration we take less than $m$ steps within the cycle to reach
the vertex $w_0$. 

We will now proceed to show that the maximum number of paths through any edge
of any of the four types is $O(m^2 2^m)$.  

\subsection{Long edges}

There are altogether $n^2$ paths, and 
for each coordinate direction $i$, half of the paths take a long edge in direction $i$.

By symmetry, as the two possible values of the other coordinates are 
not distinguished by our construction of the paths, 
all $2^{m-1}$ parallel classes of long edges are used the same number of times.
Thus for each of the parallel classes of long edges, 
exactly $(n^2/2)/2^{m-1}$ paths use (exactly) one edge from this class.
Thus at most $n^2/2^m$ paths use any individual long edge.
With $n=O(m2^m)$, we get $n^2/2^m=O(m^2 2^m)$. 

\subsection{Medium and short edges}

Each path is constructed, by the algorithm described above, 
in not more than $2m$ runs of the Procedure P, where in each iteration step
we may take a long edge, then 
we may take a few steps in a cycle of short edges, 
and then take one single step on a medium edge if needed.

Extra edges are not used at all by our canonical paths.

Thus, if all $n^2$ paths are taken together, the Procedure P
is performed not more than $2m n^2$ iterations times.

Moreover, due to the cyclic symmetry of the construction,
which treats all directions equally,
and due to the symmetry between the two values in each coordinate direction,
all $m2^m$ different (and disjoint!) edge sets of 
``a cycle of short edges plus the medium edge leaving it''
is used by \emph{the same} number of paths.
(This reflects the fact that the symmetry group of the subgraph $\CG_m'$
formed by the medium and the long edges, which is a cyclically connected cube graph,
is transitive on the medium edges.)
Thus no short edge on a cycle, or medium edge leaving it,
is used by more than $(2mn^2)/(m2^m)=2n^2/2^m$ paths.
With $n=O(m2^m)$, we get $2n^2/2^m=O(m^2 2^m)$. 
 
\subsection{Wrapping things up}

For the graphs $\CG_m''$ of the polytopes $\NCC_4(m)''$,
we have obtained in Section~\ref{subsec:NC2} that they are $4$-regular on $n=(6m-12)2^m$ 
vertices.
In Section \ref{sec:paths} we have constructed canonical paths such that no edge is used more
than $2n^2/2^m$ times.
With Sinclair's Lemma from Section~\ref{subsec:Sinclair}, this implies the expansion bound
\[
\mathcal{X}(\CG_m'') 
\ge \frac{n}{2\,\phi_{\max}} 
\ge \frac{n}{2n^2/2^m}
  = \frac{2^m}{2n}
  = \frac{2^m}{2(6m-12)2^m}
  = \frac{1}{12(m-2)}. 
\]
With the lemma from Section~\ref{subsec:expansion}
this yields that all separators of $\CG_m''$ have size at least 
\[
\frac c4 \mathcal{X}(\CG_m'')n 
\ge \frac{c\,n}{48(m-2)} 
 =  \Omega\Big(\frac n{\log n}\Big),
\]
as $n = (6m-12)2^m$ and thus $m=\Theta(\log n)$.

On the other hand, it is easy to see that separators of this order of magnitude exist:
For this, for example, look at all the $\frac nm$ long edges of one particular direction
(each vertex is incident to one long edge, and there is the same number of long edges in all directions),
and take as vertex separator one end of each of these. This yields a separator of 
size $\frac nm=O(\frac n{\log n})$. \hfill \qed

\section{Final comments}

Our construction can be extended to higher dimensions by taking prisms over $\NCC_4(m)''$. 
For a fixed dimension $d$ the resulting polytope has $2^{d-4}$ times as many vertices as the original. 
There are $d-4$ new directions and every vertex has an edge in each of these directions. 
We can amend our algorithm for example by taking edges in these directions first 
and then proceeding in the usual way. 

It is still unknown if there can be graphs of simple polytopes with even larger smallest separators.
In particular, it would be interesting to know if these graphs can be true expanders, that is, 
have smallest separators with $\Omega(n)$ vertices.

\subsection*{Acknowledgement}
Thanks to Johanna Steinmeyer for the \emph{TikZ} pictures. 
\bibliographystyle{amsplain}


\begin{thebibliography}{10}

\bibitem{joswig:_neigh}
Michael Joswig and Thilo R\"orig, \emph{Neighborly cubical polytopes and
  spheres}, Israel J. Math. \textbf{159} (2007), 221--–242.

\bibitem{Z62}
Michael Joswig and G\"unter~M. Ziegler, \emph{Neighborly cubical polytopes},
  Discrete Comput. Geometry (Gr\"unbaum Festschrift: G. Kalai, V. Klee, eds.)
  \textbf{24} (2000), 325--344.

\bibitem{Kaibel}
Volker Kaibel, \emph{On the expansion of graphs of 0/1-polytopes}, in: ``The Sharpest
  Cut: The Impact of Manfred Padberg and his Work'', MOS-SIAM Series on
  Optimization, SIAM, Philadelphia PA, 2004, pp.~199--216.

\bibitem{kalai04:_polyt}
Gil Kalai, \emph{Polytope skeletons and paths}, Handbook of Discrete and
  Computational Geometry (J.~E. Goodman and J.~{O'Rourke}, eds.), CRC Press,
  Boca Raton, second ed., 2004, First edition 1997, pp.~455--476.

\bibitem{LiTa}
Richard~J. Lipton and Robert~E. Tarjan, \emph{A separator theorem for planar
  graphs}, SIAM J. Applied Math. \textbf{36} (1979), 177--189.

\bibitem{Z150}
Lauri Loiskekoski and G\"unter~M. Ziegler, \emph{Simple polytopes without small
  separators}, Preprint, October 2015, 7~pages,
  \href{http://arxiv.org/abs/1510.00511}{arXiv:1510.00511}; {Israel J. Math.},
  to appear.

\bibitem{MillerTengThurstonVavasis}
Gary~L. Miller, Shang-Hua Teng, William Thurston, and Stephen~A. Vavasis,
  \emph{Separators for sphere-packings and nearest neighbor graphs}, J. ACM
  \textbf{44} (1997), 1--29.

\bibitem{PreparataVuillemin1981}
Franco~P. Preparata and Jean Vuillemin, \emph{The cube-connected cycles: a
  versatile network for parallel computation}, Communications of the ACM
  \textbf{24} (1981), no.~5, 300--309.

\bibitem{ringel55:_ueber_probl_wuerf_wuerf}
Gerhard Ringel, \emph{{\"Uber drei kombinatorische Probleme am
  $n$-dimensionalen W\"urfel und W\"urfel\-gitter}}, Abh. Math. Sem. Univ.
  Hamburg \textbf{20} (1955), 10--19.

\bibitem{Z102}
Raman Sanyal and G\"unter~M. Ziegler, \emph{Construction and analysis of
  projected deformed products}, Discrete Comput. Geometry \textbf{43} (2010),
  412--435.

\bibitem{Sinclair}
Alistair Sinclair, \emph{{Algorithms for Random Generation and Counting: A
  Markov Chain Approach}}, Progress in Theoretical Computer Science,
  Birkhäuser, Boston, 1993.

\bibitem{Z100}
G\"unter~M. Ziegler, \emph{Polyhedral surfaces of high genus}, in: ``Discrete
  Differential Geometry'', Oberwolfach Seminars, vol.~38, Birkh\"auser, Basel
  2008, pp.~191--213.

\end{thebibliography}


\end{document}